\documentclass[12pt]{amsart}
\usepackage{amsthm,amsmath,amssymb,amscd,graphicx,enumerate, stmaryrd,xspace,verbatim, epic, eepic,color,url,caption,subcaption}

\usepackage[active]{srcltx}
\usepackage[all]{xypic}
\SelectTips{cm}{}

{
   \newtheorem{theorem}{Theorem}
      \newtheorem*{theorem*}{Theorem}
   \newtheorem{proposition}{Proposition} 

   \newtheorem*{conjecture*}{Conjecture}
   
}
{\theoremstyle{definition}
          \newtheorem*{exercise*}{Exercise}

   \newtheorem*{example*}{Example}

   \newtheorem*{definition*}{Definition}

}
\newcommand{\RR}{{\mathbb{R}}}

\newcommand{\CC}{{\mathbb{C}}}

\newcommand{\PP}{{\mathbb{P}}}
\newcommand{\ZZ}{{\mathbb{Z}}}

\renewcommand{\AA}{{\mathbb{A}}}

\newcommand{\cM}{{\mathcal M}}

\def\<{\langle}
\def\>{\rangle}

\newcommand{\Spec}{\operatorname{Spec}}

\newcommand{\Aut}{{\operatorname{Aut}}}

\newcommand{\ocM}{\overline{{\mathcal M}}}

\newcommand{\val}{{\operatorname{val}}}

\newcommand{\Trop}{{\operatorname{Trop}}}

\newcommand{\An}{{\operatorname{An}}}

\newcommand{\double}{\genfrac..{0pt}1
{\raise -2pt\hbox{$\scriptstyle\longrightarrow$}}{\raise 4pt\hbox
{$\scriptstyle\longrightarrow$}}} 

\renewcommand{\setminus}{\smallsetminus}

\def\tototi{\mathbin{\mathop{\otimes}\limits^{\raise-1pt\hbox
{$\scriptscriptstyle {\rm L}$}}}}

\def\indlim{\mathop{\vrule width0pt height7pt depth
4pt\smash{\lim\limits_{\raise 1pt\hbox to 14.5pt
{\rightarrowfill}}}}}
\def\projlim{\mathop{\vrule width0pt height7pt depth
4pt\smash{\lim\limits_{\raise 1pt\hbox to 14.5pt
{\leftarrowfill}}}}}

\newcommand\displaceamount{3pt}

\newcommand{\doubledown}{\ar@<\displaceamount>[d]\ar@<-\displaceamount>[d]}

\newcommand{\doubleup}{\ar@<\displaceamount>[u]\ar@<-\displaceamount>[u]}

\newcommand{\doubleright}{\ar@<\displaceamount>[r]\ar@<-\displaceamount>[r]}

\begin{document}

\title{Moduli  of algebraic and tropical curves}

\author[Abramovich]{Dan Abramovich}



\address{Department of Mathematics\\
Brown University\\
Box 1917\\
Providence, RI 02912\\
U.S.A.}
\email{abrmovic@math.brown.edu}



\thanks{Research of Abramovich supported in part by NSF grant DMS-0901278. }
\date{\today}

\begin{abstract}
 {\it Moduli spaces} are a geometer's obsession.
A celebrated example in algebraic geometry is the space $\bar M_{g,n}$ of stable
$n$-pointed algebraic curves of genus $g$, due to Deligne--Mumford and Knudsen. It has a
delightful combinatorial structure based on {\em weighted graphs}.

Recent papers of Branetti, Melo, Viviani and of Caporaso defined an entirely different
moduli space of {\em tropical curves}, which are weighted  metrized graphs. It also
has a delightful combinatorial structure based on weighted graphs.

One is led to ask whether there is a geometric connection between these moduli spaces.
In joint work  \cite{ACP} with Caporaso and Payne, we exhibit a connection, which
passes through a third type of geometry - nonarchimedean analytic geometry.

\end{abstract}
\maketitle
\section{The moduli bug}

Geometers of all kinds are excited, one may say obsessed, with moduli spaces; these are the spaces which serve as parameter spaces for the basic spaces geometers are most interested in. 

It was none other than Riemann who introduced the moduli bug into geometry, when he noted that Riemann surfaces of genus $g>1$ ``depend on $3g-3$ {\em Moduln}". This is a consequence of his famous ``Riemann existence theorem", which tells us how to put together a Riemann surfaces by slitting a number of copies of the Riemann sphere and sewing them together. The number $3g-3$ is simply the number of ``effective complex parameters" necessary for obtaining {\em every} Riemann surface this way. \footnote{It was none other than Riemann who at the same stroke introduced the word ``moduli", with which we have been stuck ever since.}

This brings us to the classic case of the moduli phenomenon, and the obsession that comes with it, the space $\cM_g$ of Riemann surfaces of genus $g$: fixing a compact oriented surface $S$, each point on $\cM_g$ corresponds uniquely to a complex structure $C$ on $S$. Being an algebraic geometer, I tend to think about these Riemann surfaces as ``smooth projective and connected complex algebraic curves of genus $g$", or just ``curves of genus $g$" in short. This explains the choice of the letter $C$.\footnote{I hope I can be excused for the confusion between ``surfaces" and ``curves", which comes from the fact that the  dimension of $\CC$ as a real manifold is 2.}

 The moduli space  $\cM_g$ is the result of important work of  many mathematicians, such as Ahlfors, Teichm\"uller, Bers, Mumford ...

\section{From elliptic curves to higher genus} 

The first example, of elliptic curves, is familiar from Complex Analysis, where an elliptic curve is defined as the quotient $\CC/\langle 1, \tau\rangle$ of the complex plane by a lattice of rank 2 with $Im(\tau)>0$. We learn, for instance in Ahlfors's book \cite{Ahlfors}, that isomorphism classes of elliptic curves are identified uniquely by the so called $j$-invariant $j(\tau)$, an important but complicated analytic function on the upper half plane. In algebraic geometry one can use a shortcut to circumvent this: every elliptic curve has a so called Weierstrass equation $$E_{a,b}: \ y^2 = x^3 + ax + b$$ with nonzero discriminant $\Delta(a,b) = 4a^3 + 27b^2\neq 0$. One can identify the $j$-invariant as $$j(a,b) = \frac{4a^3}{4a^3 + 27b^2} \in \CC,$$ so that two elliptic curves are isomorphic: $E_{a,b} \simeq E_{a',b'}$ if and only if $j(a,b) = j(a',b')$.

\begin{figure}[htb]
\begin{center}
\includegraphics[width=8cm]{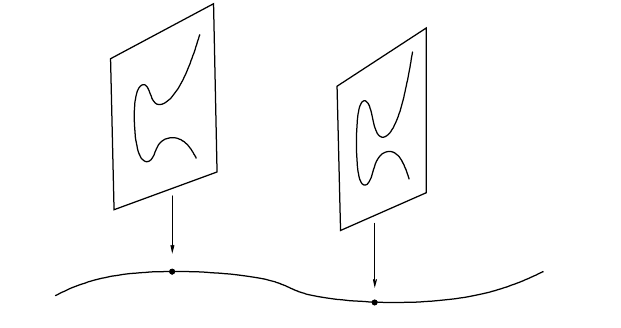}
\caption{The family of elliptic curves over $\CC$}\label{Fig:j-line}
\end{center}
\end{figure}

Either way, the moduli space of elliptic curves is just $\CC$ - or, in the language of algebraic geometers, the affine line $\AA^1_\CC$ (Figure \ref{Fig:j-line}). 

The story for genus $g>1$ is quite a bit more involved. But the principle, at least in algebraic geometry, is similar: note that $j(a,b)$ is an invariant rational function in the parameters $a,b$, which are the coefficients of the defining equation of $E_{a,b}$ written in its Weierstrass form. For higher genus one does the same: one finds a sort of canonical form for a Riemann surface in a suitable projective space, one collects the coefficients of the defining equations,  and then the coordinates on $M_g$ are invariant rational functions in these.  The result, in its algebraic version due to Mumford, is:

\begin{theorem} The space $M_g$ is a complex quasi projective variety.
\end{theorem}

It is a rather nice variety - it is not quite a manifold, but it is an {\em orbifold}: it locally looks like the quotient of a manifold by the action of a finite group.

In general the global geometry of $\cM_g$ is quite a bit more involved than the geometry of $\CC$. Its complex dimension is indeed $3g-3$.

\section{The problem of compactness}

Angelo Vistoli from Pisa has said that ``working with a noncompact space is like trying to keep your change when you have holes in your pockets". The space $\CC$ of elliptic curve, and the space $\cM_g$ of curves of genus $g$, are noncompact, and one wishes to find a natural compactification. 

Of course every quasi-projective variety sits inside a projective space, and its closure is a compactification. But that is not natural: we want a compactification which is itself  a moduli space, of slightly singular Riemann surfaces!

For instance, the moduli space of elliptic curves $\CC$ has a nice compactification $\PP^1_\CC$, the Riemann sphere. In which way does the added point $\infty$ represent a singular Riemann surface?

\begin{figure}[htb]
\begin{center}
\includegraphics[width=8cm]{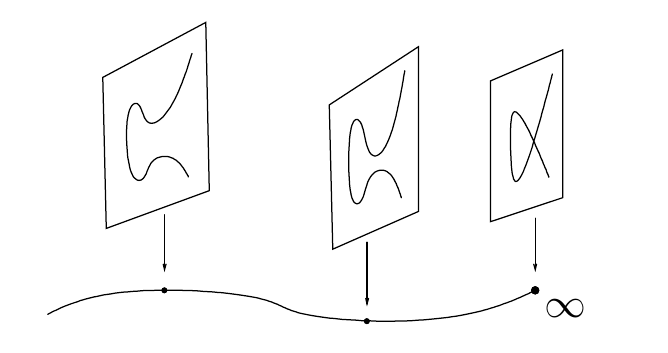}
\caption{The family of elliptic curves over $\PP^1_\CC$}\label{j-line-cpt}
\end{center}
\end{figure}

The function $j(a,b)$ extends to a regular function $j: \CC^2\setminus \{(0,0)\} \to \PP^1$, and indeed all the points where $\Delta = 0$ in   $\CC^2\setminus \{(0,0)\}$ correspond to ``singular elliptic curves" such as $E_{-3,2}:  y^2 = x^3 -3x + 2$: the singular point $(x,y) = (1,0)$ has local plane coordinates with equation of the form $zw=0$, a so called {\em node}.  In fact all these singular elliptic curves are isomorphic! (See Figure \ref{j-line-cpt}.)

\begin{figure}[htb]
\begin{center}

\includegraphics[width=8cm]{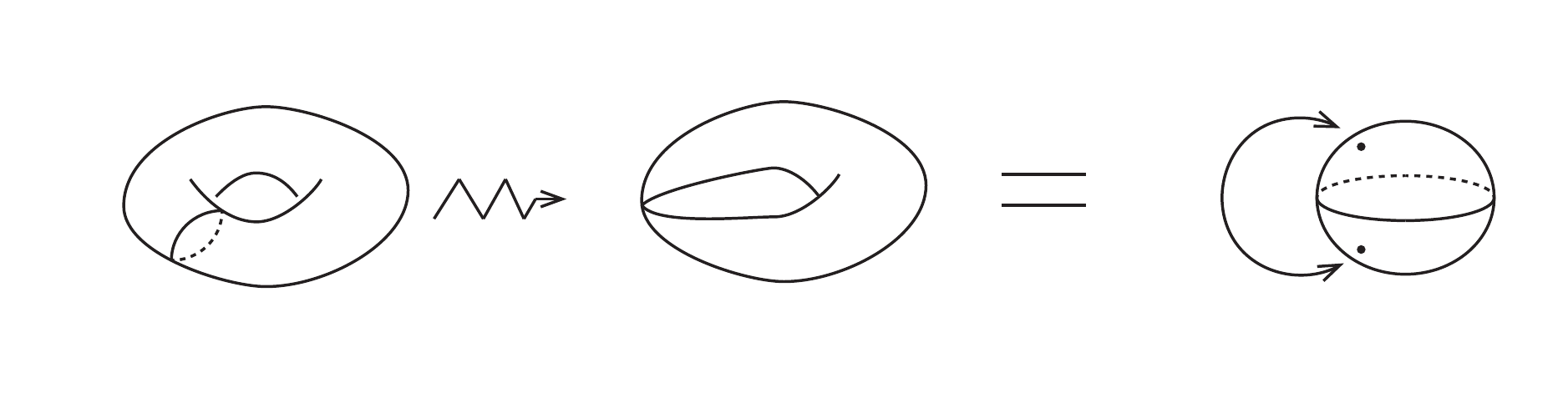}
\caption{A degenerate elliptic curve as a sphere with glued points }\label{deg-eg}
\end{center}
\end{figure}

From the Riemann surface point of view one can describe such singular Riemann surfaces as follows: choose your favorite elliptic curve $E$, thought of as a torus, and wrap a loop around its girth. Now pull the loop until it pops. What you get is the Riemann sphere with two points glued  together (Figure \ref{deg-eg}).

\begin{figure}[htb]
\begin{center}

\includegraphics[width=8cm]{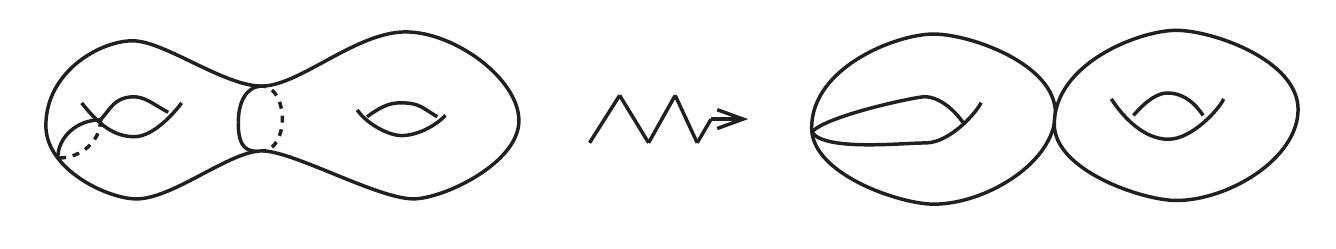}
\caption{A degenerate Riemann surface of genus 2}\label{deg-g2}
\end{center}
\end{figure}

Deligne and Mumford showed that this can be done in general: they described degenerate algebraic curves of genus $g$, obtained by choosing a number of disjoint loops and pulling them until they pop (Figure \ref{deg-g2}).

\begin{figure}[htb]
\begin{center}
\includegraphics[width=6cm]{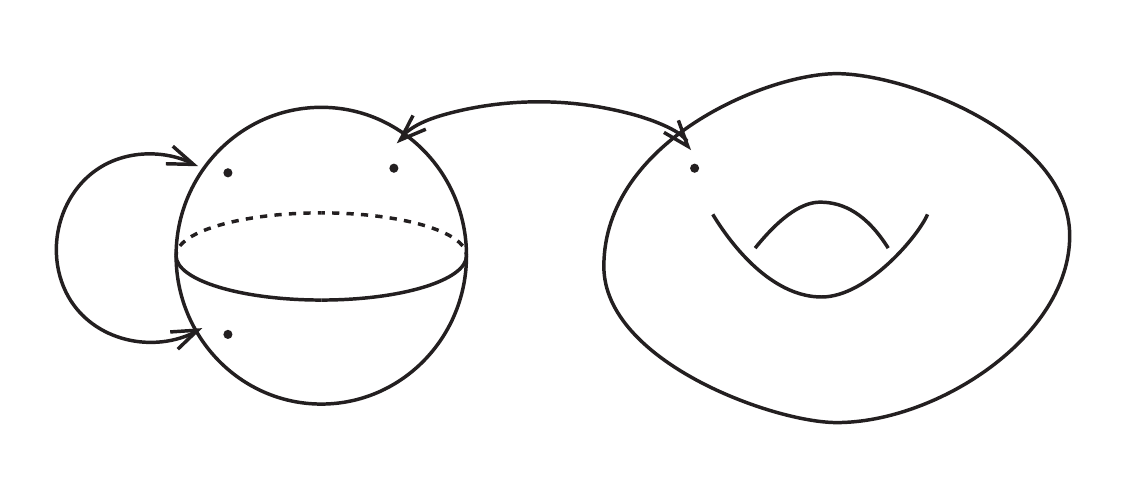}
\caption{Gluing the same degenerate Riemann surface of genus 2}\label{glue-g2}
\end{center}
\end{figure}

 The result is a singular Riemann surface obtained by taking a number of usual Riemann surfaces, choosing a number of points of them, and indicating how these points are to be glued together (Figure \ref{glue-g2}).
 
  One imposes a {\em stability condition} on the loops chosen, which is best described in combinatorial terms, see below.
 The collection of all of these objects  is Deligne and Mumford's moduli space of {\em stable curves} $\ocM_g$. It is also a {\em projective} orbifold containing $\cM_g$ as a dense open subset.

\section{The weighted graph of a curve}

The combinatorial underpinning of a stable curve is given by its {\em dual graph} $\Gamma$. This is a so called {\em weighted graph} where each vertex $v$ is assigned an integer weight $g(v) \geq 0$.

\begin{figure}[htb]
\begin{center}
\begin{subfigure}[c]{2in}\includegraphics[width=4cm]{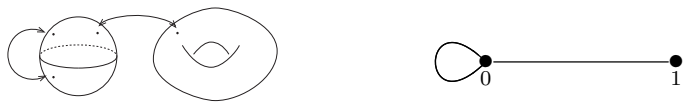}\end{subfigure}
\hspace{1cm}
\begin{subfigure}[b]{2in}
\xymatrix@=5.5pc
{*{\bullet}\ar@{{-}{-}{-}}@(ul,dl)[]\ar@{-}[r]_(1)1_(0)0&*{\bullet}}\end{subfigure}
\caption{The glued curve ... and its graph}\label{graph-g2}
\end{center}
\end{figure}

 Given a singular Riemann surface $C$ as above, its graph has a vertex $v_i$ corresponding to each component $C_i$, weighted by the genus $g(C_i)$. Corresponding to each node, where a point of $C_i$ is glued to a point of $C_j$ one writes an edge tying $v_i$ to $v_j$ (Figure \ref{graph-g2}). 
 
 The genus of $\Gamma$, and of any corresponding singular Riemann surface, is given by a simple formula involving the first Betti number of the graph: $$b_1(\Gamma,\ZZ) + \sum_{v\in V(\Gamma)} g(v).$$ 
 
 We can now describe the stability condition: the graph $\Gamma$, and any corresponding curve, is {\em stable} if every vertex $v$ of genus $g(v) = 0$ has valence ${\val}(v) \geq 3$ and every vertex $v$ of genus $g(v) = 1$ has valence $\val(v) \geq 1$.
 
 \section{The combinatorial structure of moduli space}
 
These graphs give us a way to put together the space $\ocM_g$ piece by piece. 

For each weighted graph $\Gamma$ there is a nice moduli space $\cM_\Gamma$ parametrizing singular Riemann surfaces with  weighted graph $\Gamma$.  Each $\cM_\Gamma$ is an orbifold, and its codimension in $\ocM_g$ is simply the number of edges $|E(\Gamma)|$. We have $$\ocM = \mathop{\coprod}\limits_{g(\Gamma)=g} \cM_\Gamma.$$

\begin{figure}[htb]
\begin{center}
\includegraphics[width=8cm]{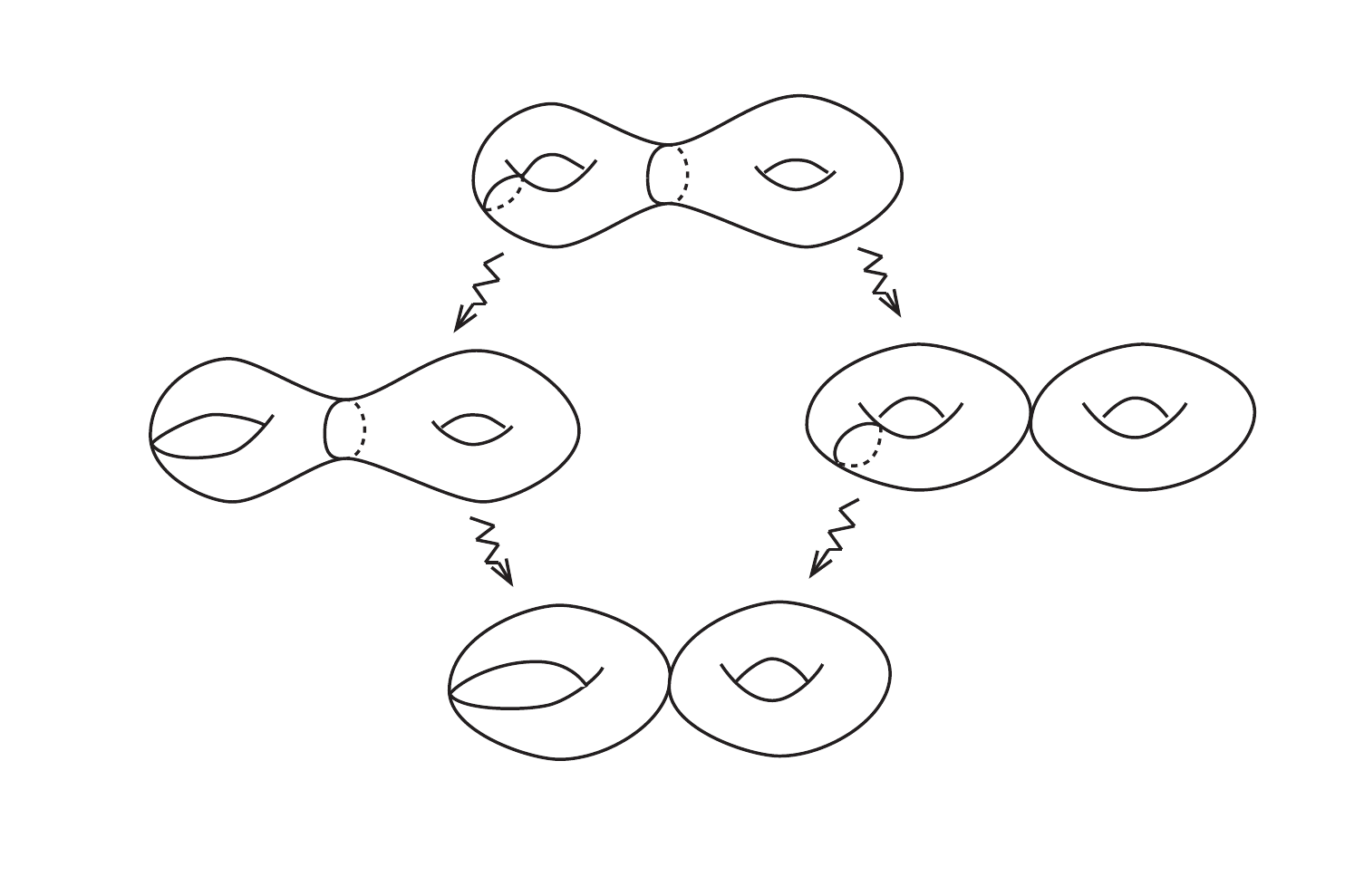}
\caption{Step-by-step degenerations in genus 2}\label{degs-g2}
\end{center}
\end{figure}

The pieces $\cM_\Gamma$  form a {\em stratification} of $\cM_g$, in the sense that the closure of $\cM_\Gamma$ is the disjoint union of pieces of the same kind. To determine the combinatorial structure of $\ocM_g$ we need to record which pieces $\cM_{\Gamma'}$ appear in the closure: these correspond to singular Riemann surfaces $C'$ where {\em more} loops were pulled until they popped than in a curve $C$ corresponding to $\Gamma$ (Figure \ref{degs-g2}).

\begin{figure}[htb]
\begin{center}
\xymatrix@=0.5pc{&&&&&&&&&&&&&&&&&&&&&&&\\
&&*{\bullet}\ar@{{-}{.}}[llu]\ar@{{-}{.}}[ll]\ar@{{-}{.}}[lld]\ar@{-}[r]_(0){g_1}_(1){g_2} & *{\bullet}\ar@{{-}{.}}[rru]\ar@{{-}{.}}[rr]\ar@{{-}{.}}[rrd]&&\ar@{~>}[rr]&&&&
*{\bullet}\ar@{{-}{.}}[llu]\ar@{{-}{.}}[ll]\ar@{{-}{.}}[lld] \ar@{{-}{.}}[rru]\ar@{{-}{.}}[rr]\ar@{{-}{.}}[rrd]_(0.3){g_1+g_2}&&&&&
&&*{\bullet}\ar@{{-}{-}{-}}@(ul,dl)[] \ar@{{-}{.}}[rru]\ar@{{-}{.}}[rr]\ar@{{-}{.}}[rrd]_(0.1){g}&&\ar@{~>}[rr]&&&
*{\bullet}  \ar@{{-}{.}}[rru]\ar@{{-}{.}}[rr]\ar@{{-}{.}}[rrd]_(0.3){g+1}&&
\\
&&&&&&&&&&&&&&&&&&&&&&&}
\caption{Contracting an edge ... and a loop}\label{contr}
\end{center}
\end{figure}
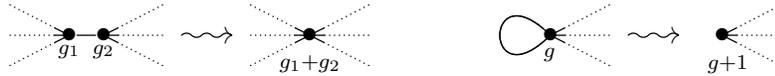

On the level of graphs this corresponds to simply saying that there is a {\em contraction} $\Gamma' \to \Gamma$: on weighted graphs, contracting an edge connecting vertices $v_1,v_2$ with genera $g_1,g_2$ results in a vertex with genus $g_1+g_2$; similarly contracting a loop on a vertex with genus $g$ results with a vertex with genus $g+1$ (Figure \ref{contr}).

 So the combinatorial structure of $\ocM_g$ is given by the following rule:

$$\cM_{\Gamma'} \subset \ocM_\Gamma \ \ \ \Longleftrightarrow \ \ \ \exists \text{ contraction } \ \Gamma' \to \Gamma.$$

\begin{figure}[htb]
\begin{center}
\begin{subfigure}[c]{2.5in}
\includegraphics[width=5.5cm]{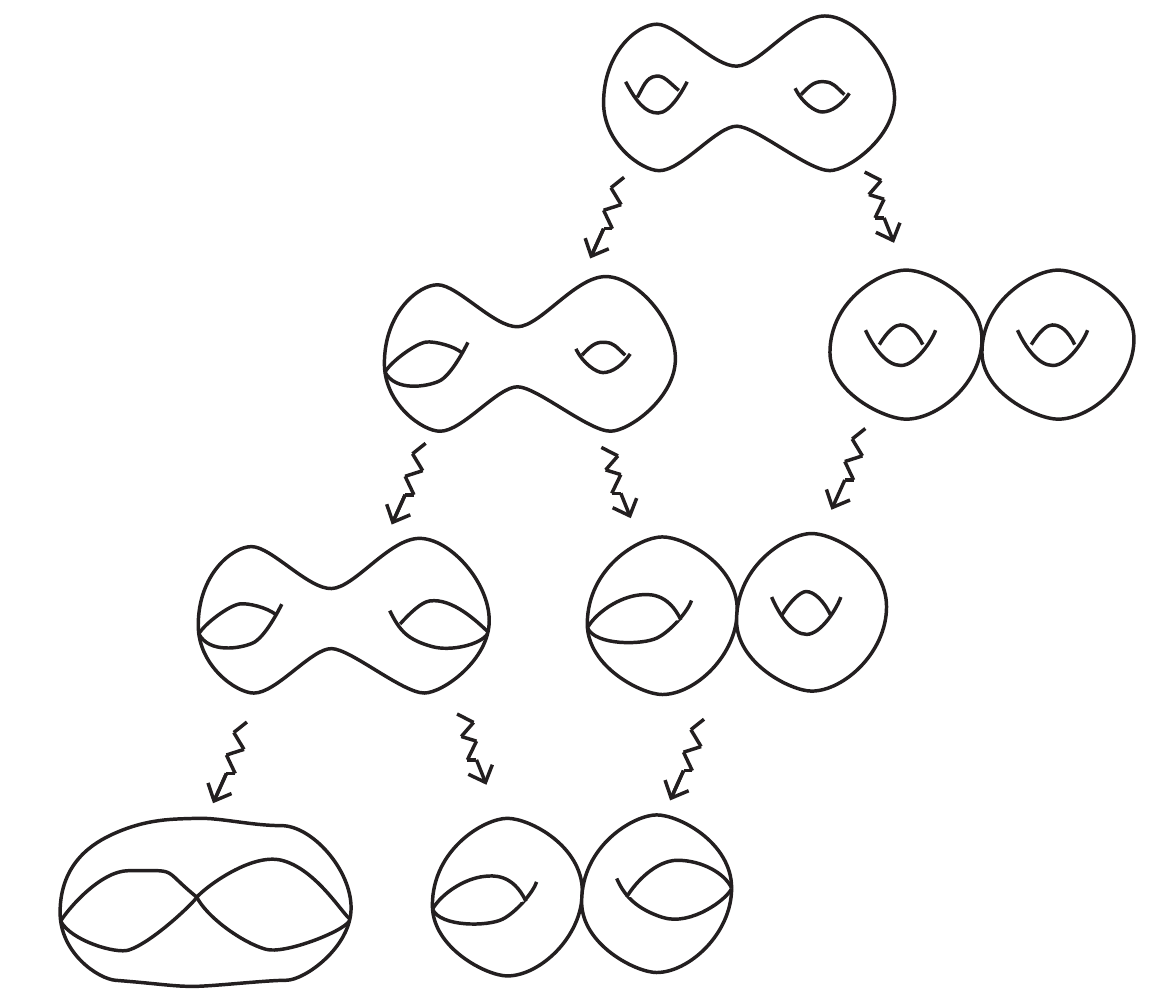}
\end{subfigure}
\begin{subfigure}[c]{3in}
\xymatrix@=0.4pc{\\ \\
&&&&&&&&*{\bullet}\ar@{}[r]_(0)2&&&&\\
&&&&&&&&&&&&\\
&&&&&&\ar@{~>}[ru]&&&&\ar@{~>}[lu]\\
&&&&&&*{\bullet}\ar@{{-}{-}{-}}@(lu,ld)[]\ar@{}[r]_(0)1&&&&*{\bullet}\ar@{{-}{-}{-}}[r]_(0)1_(1)1&*{\bullet}\\
&&&&&&&&&&&&\\
&&&&\ar@{~>}[ru]&&&\ar@{~>}[lu]&&\ar@{~>}[ru]\\
&&&*{\bullet}\ar@{{-}{-}{-}}@(lu,ld)[]\ar@{{-}{-}{-}}@(ru,rd)[]\ar@{}[r]_(0)0&&&&&*{\bullet}\ar@{{-}{-}{-}}@(lu,ld)[]\ar@{{-}{-}{-}}[r]_(0)0_(1)1&*{\bullet}\\
&&&&&&&&&&&&\\
*{\bullet}\ar@{{-}{-}{-}}@/_10pt/[dd]\ar@{{-}{-}{-}}@/^10pt/[dd]\ar@{{-}{-}{-}}[dd]\ar@{}[r]^(0)0&\ar@{~>}[ru]&&&&\ar@{~>}[lu]&&\ar@{~>}[ru]\\
&&&&&&*{\bullet}\ar@{{-}{-}{-}}@(lu,ld)[]\ar@{{-}{-}{-}}[r]_(0)0_(1)0&*{\bullet}\ar@{{-}{-}{-}}@(ru,rd)[]\\
*{\bullet}\ar@{}[r]_(0)0&
}
\end{subfigure}
\caption{Curves in $\ocM_2$ ... and their graphs}\label{Fig:m2}
\end{center}
\end{figure}

The skeletal picture of $\ocM_2$ is given in Figure \ref{Fig:m2}. The top line is the big stratum $\cM_2$ of complex dimension 3, and the bottom strata are points, of dimension $0$.

\section{Tropical curves}
There is another geometry which builds on the combinatorics of weighted graphs, namely the geometry of {\em tropical curves}. This is a much more recent theory. One can identify its pre-history in the work of Culler--Vogtmann on outer space \cite{Culler-Vogtmann}, and continuing with the work of Mikhalkin on tropical enumeration of plane curves \cite{Mikhalkin,Mikhalkin1}. The notion of tropical curves in the sense described here was introduced by Brannetti--Melo--Viviani \cite{BMV}  and Caporaso \cite{C2}. \footnote{The name ``tropical" is the result of tradition: tropical curves arise in tropical geometry, which sometimes relies on min-plus algebra. This was studied by the Brazilian (thus tropical) computer scientist Imre Simon.}

A {\em tropical curve} is simply a {\em metric} weighted graph $$G = (\Gamma, \ell:E(\Gamma) \to \RR_{>0}\cup\{\infty\}).$$ It consists of a weighted graph $\Gamma$ and a possibly infinite length $\ell(e)>0$  assigned to each edge. 

\section{Moduli of tropical curves}

\begin{figure}[htb]
\begin{center}
\xymatrix@=0.5pc{&&&&&&&&&&&&&&&&&&&&&&&&&&\\
&&*{\bullet}\ar@{{-}{.}}[llu]\ar@{{-}{.}}[ll]\ar@{{-}{.}}[lld]\ar@{-}[r]_(0){g_1}_(1){g_2}^t & *{\bullet}\ar@{{-}{.}}[rru]\ar@{{-}{.}}[rr]\ar@{{-}{.}}[rrd]&&\ar@{->}[rrr]_{t\to 0}&&&&&
*{\bullet}\ar@{{-}{.}}[llu]\ar@{{-}{.}}[ll]\ar@{{-}{.}}[lld] \ar@{{-}{.}}[rru]\ar@{{-}{.}}[rr]\ar@{{-}{.}}[rrd]_(0.3){g_1+g_2}&&&&&
&&*{\bullet}\ar@{{-}{-}{-}}@(ul,dl)[]_t \ar@{{-}{.}}[rru]\ar@{{-}{.}}[rr]\ar@{{-}{.}}[rrd]_(0.1){g}&&\ar@{->}[rrr]_{t\to 0}&&&&
*{\bullet}  \ar@{{-}{.}}[rru]\ar@{{-}{.}}[rr]\ar@{{-}{.}}[rrd]_(0.3){g+1}&&
\\
&&&&&&&&&&&&&&&&&&&&&&&&&&}
\caption{Pulling an edge ... and a loop}\label{limit}
\end{center}
\end{figure}
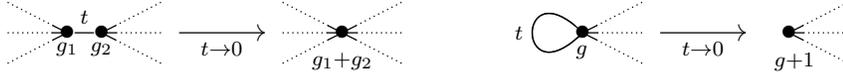

Being a geometer, one is infected with the moduli bug. Therefore the moment one writes the definition of a tropical curve one realizes that they have a moduli space. Fixing a weighted graph $\Gamma$, the tropical curves having graph isomorphic to $\Gamma$ are determined by the lengths of the edges, and the collection of lengths is unique up to the permutations obtained by automorphisms of the graph. We can therefore declare the moduli space of such tropical curves to be
$$\cM_\Gamma^\Trop\ \ \ \  = \ \ \ \ (\RR_{>0}\cup\{\infty\})^{|E(\Gamma)|}\ \ /\ \ \Aut(\Gamma).$$
We can put together these moduli spaces $\cM_\Gamma^\Trop$ by observing the following: if we let the length of an edge $e$ in $G$ approach $0$, the metric space $G$ approaches $G'$, which is the metric graph associated to the graph $\Gamma'$ obtained by contracting $e$, as in Figure \ref{limit}.

This allows us to glue together $\cM_\Gamma^\Trop$ in one big moduli space 
$$\ocM_g^\Trop = \mathop{\coprod}\limits_{g(\Gamma)=g}  \cM_\Gamma^\Trop.$$
It is a nice compact cell complex. 

Note that the gluing rule precisely means that

$$\cM^\Trop_{\Gamma'} \supset \overline{\cM^\Trop_\Gamma} \ \ \ \Longleftrightarrow \ \ \ \exists \text{ contraction } \ \Gamma' \to \Gamma.$$

\section{The question of comparison}

We obtained two geometries associated to the combinatorics of weighted graphs, summarized as follows:

$$\xymatrix{
\ocM_g\ar@{~>}[dr] && \ocM_g^\Trop\ar@{~>}[dl] \\ &\{\Gamma\}
}$$

The moduli spaces are clearly similar in their combinatorial nature, and one would like to tie them together somehow:

$$\xymatrix{
\ocM_g\ar@{~>}[dr]\ar@{.}@/^2pc/[rr]|-{?} && \ocM_g^\Trop\ar@{~>}[dl] \\ &\{\Gamma\}
}$$

The question was raised in \cite[Section 6]{C2}. A possible answer was suggested by Baker-Payne-Rabinoff \cite[Remark 5.52]{BPR}, Tyomkin \cite[Section 2]{Tyomkin}, Viviani \cite[Theorem A]{Viviani}. Below I report on joint work \cite{ACP} with Caporaso and Payne, where we prove the suggested answer to be correct; it is based on the non-archimedean analytic spaces  of Berkovich \cite{Berkovich} and their skeleta, and fundamental constructions of such skeleta by Thuillier \cite{Thuillier}.

\begin{figure}[htb]
\begin{center}
$$\xymatrix@=0.4pc{\text{dimension}\\
0\ar@{.}[rrrrrrrrr]&&&&&&&&&&*{\bullet}\ar@{}[r]_(0)2&&&&\\
&&&&&&&&&&&&\\
&&&&&&&&\ar@{~>}[ru]&&&&\ar@{~>}[lu]\\
1\ar@{.}[rrrrrr]&&&&&&&&*{\bullet}\ar@{{-}{-}{-}}@(lu,ld)[]\ar@{}[r]_(0)1&&&&*{\bullet}\ar@{{-}{-}{-}}[r]_(0)1_(1)1&*{\bullet}\\
&&&&&&&&&&&&&\\
&&&&&&\ar@{~>}[ru]&&&\ar@{~>}[lu]&&\ar@{~>}[ru]\\
2\ar@{.}[rrr]&&&&&*{\bullet}\ar@{{-}{-}{-}}@(lu,ld)[]\ar@{{-}{-}{-}}@(ru,rd)[]\ar@{}[r]_(0)0&&&&&*{\bullet}\ar@{{-}{-}{-}}@(lu,ld)[]\ar@{{-}{-}{-}}[r]_(0)0_(1)1&*{\bullet}\\
&&&&&&&&&&&&&\\
&&*{\bullet}\ar@{{-}{-}{-}}@/_10pt/[dd]\ar@{{-}{-}{-}}@/^10pt/[dd]\ar@{{-}{-}{-}}[dd]\ar@{}[r]^(0)0&\ar@{~>}[ru]&&&&\ar@{~>}[lu]&&\ar@{~>}[ru]\\
3\ar@{.}[r]&&&&&&&&*{\bullet}\ar@{{-}{-}{-}}@(lu,ld)[]\ar@{{-}{-}{-}}[r]_(0)0_(1)0&*{\bullet}\ar@{{-}{-}{-}}@(ru,rd)[]\\
&&*{\bullet}\ar@{}[r]_(0)0&
}$$
\caption{Graph contractions in genus 2}\label{Fig:m2-trop}
\end{center}
\end{figure}
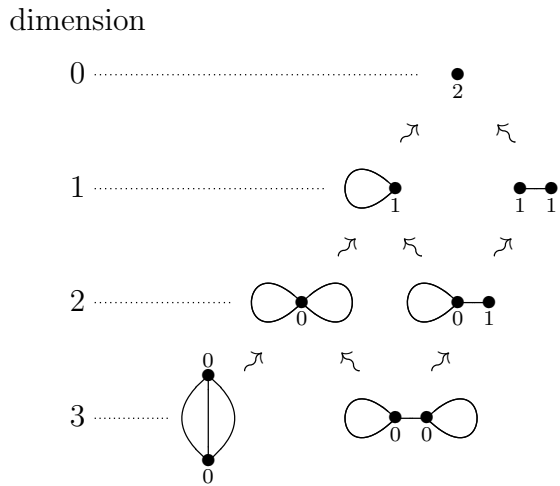

One slightly disturbing feature is the fact that the combinatorial structures - the stratifications - of the moduli spaces $\ocM_g$ and $\ocM_g^\Trop$ are reversed!  In Figure \ref{Fig:m2-trop}, the top startum is a point in dimension 0, and the bottom strata are of dimension 3.

\section{Nonarchimedean analytic geometry} 

A {\em valued field} is a field $K$ with a multiplicative seminorm $K \to \RR_{\geq 0}$. One often translates the seminorm to a {\em valuation} $v: K \to \RR\cup\{\infty\}$ by declaring $\val(x) = -\log \| x\|$. The valued field is {\em nonarchimedean} if, like the field of $p$-adic numbers, it satisfies the strict triangle inequality $$\| a+b\|  \leq \max(\| a\| ,\| b\| ).$$ One defines the valuation ring to be $$R = \{x\in K: \| x\| \leq 1\} = \{x\in K: \val(x) \geq 0\},$$ which has a unique maximal ideal $I = \{x\in K: \| x\| < 1\}$. The residue field is defined as $\kappa = R/I$.

Recall that in scheme theory, a point of a variety $X$ corresponds to a field extension $K\supset \CC$ and a morphism $\Spec K \to X$, up to equivalence given by further extensions. 
In particular a point of $\ocM_g$ corresponds to a field extension $K\supset \CC$ and a stable curve $C/\Spec K$. Topologists cannot be happy about this structure, because the Zariski topology of a scheme is not Hausdorff in the least.

Berkovich associates to $X$ an analytic variety $X^\An$ - a locally ringed space which admits a natural morphism $X^\An \to X$. A point of 
$X^\An$ corresponds to a  nonrchimedean valued field extension $K\supset \CC$ (extending the {\em trivial} valuation on $\CC$) and a morphism $\Spec K \to X$, up to equivalence by further valued field extensions.

In particular a point of $\ocM_g^\An$ corresponds to a {\em valued} field extension $K\supset \CC$ and a stable curve $C/\Spec K$. 

Since every valued field extension is, in particular, a field extension, there is a morphism $\ocM_g^\An \to \ocM_g$.  It is a bit of magic that adding the valuations ``stretches'' the generic points just enough to make the space Hausdorff and locally connected. 
 We can now extend our diagram of relationships as follows:

$$\xymatrix{ 
& \ocM_g^\An \ar[dl] \ar@{.}[dr]|-{?} \\
\ocM_g\ar@{~>}[dr] && \ocM_g^\Trop\ar@{~>}[dl] \\ &\{\Gamma\}
}$$

\section{Making the connection}

We are finally ready to close the diagram. 

We said that  a point of $\ocM_g^\An$ corresponds to a  valued field extension $K\supset \CC$ and a stable curve $C/\Spec K$. Since $\ocM_g$ is proper, a stable curve over $K$ uniquely extends to $C / \Spec R$ over the valuation ring (at least after a further field extension, which we may ignore because of our equivalence relation). The fiber of $C / \Spec R$ over the residue field is a stable curve $C_s$ with dual graph $\Gamma(C_s)$. We need to put a metric on this graph.

Consider an edge $s\in E(\Gamma(C_s))$ corresponding to a node $p\in C_s$. Near $p$ the curve $C$ admits a local equation of the shape $xy = f$, where $f\in R$. Define $\ell(e) = \val(f)$, which is independent of the local equation chosen (or the field $K$). This results with a tropical curve $G := (\Gamma,\ell)$.

\begin{theorem}
The resulting map $\Trop: \ocM_g^\An \to \ocM_g^\Trop$ is proper, continuous and surjective.  It also makes $\ocM_g^\Trop$ canonically isomorphic to the Berkovich skeleton of $\ocM_g^\An$.
\end{theorem}

We obtain a picture as follows:

$$\xymatrix{ 
& \ocM_g^\An \ar[dl] \ar[dr] \\
\ocM_g\ar@{~>}[dr] && \ocM_g^\Trop\ar@{~>}[dl] \\ &\{\Gamma\}
}$$

In other words, we tied together the algebraic geometry of the space $\ocM_g$  of stable curves with the metric geometry of the space $\ocM_g^\Trop$ of tropical curves, by going through nonarchimedean analytic geometry. What enabled us to relate the spaces  $\ocM_g$  and $\ocM_g^\Trop$ with ``reversed" combinatorial structures was the use of the {\em reduction} of a curve over $K$ to a curve over the residue field.

\section{Comments on the proof}

A key tool in our proof is a result of Thuillier \cite{Thuillier}.  To a complex toroidal embedding $X$, with or without self intersection, Thuillier assigns a so called {\em compactified fan} $\Sigma_X$, generalizing the cone complex construction of \cite{KKMS}. This is a canonical construction of a Berkovich skeleton in the special case of toroidal varieties over a field with trivial valuation (such as $\CC$). The combinatorial structure of the fan mirrors that of the toroidal structure: there is a cone $\sigma_F\in \Sigma_X$ assigned to each stratum $F\in X$, and these are glued to each other via a natural rule which in particular says that $\sigma_{F}$ is glued as a face of $\sigma_{F'}$ precisely if $F' \subset \overline F$ - note the reversal of order!

Furtheremore, considering the associated restricted Berkovich space $X^\beth$, one has that $\Sigma_X$ canonically sits inside $X^\beth$ as a deformation retract, in particular there is a canonical proper continuous map $p:X^\beth \to \Sigma_X$. This map is described via reduction: if $K$ is a valued field and $x:\Spec K \to X$ is a point on $X^\beth$, consider the reduction $\bar x : \Spec \kappa \to X$, and assume $\bar x$ lands in the toroidal stratum $F$. Then $p(x)$ lands in  the cone $\sigma_F$ corresponding to the stratum $F$. The position of $p(x)$ on $\sigma_F$ is determined by by applying the toroidal valuations at $F$ to $x$.

This does not quite apply to the moduli space, but almost. The point is that the moduli space is not toroidal, but the moduli stack is. The technical point we needed to prove is the following slight generalization of Thuillier's result: 

\begin{proposition}
Suppose $X$ is the coarse moduli space of a toroidal stack. Then there is a canonical compactified fan $\Sigma_X$ which is a deformation retract of $X^\beth$. The map $p: X^\beth \to \Sigma_X$ admits the same description in terms of reductions as above.
\end{proposition}

When $X$ is proper we have $X^\beth=X^\An$. 
This applies to the moduli space $\ocM_g$, so $\ocM_g^\beth = \ocM_g^\An$.

We obtain the following picture

$$\xymatrix{ 
 \ocM_g^\An \ar^p[r] \ar_{\Trop}[dr] & \Sigma_X\ar@{.}|-{?}[d] \\
& \ocM_g^\Trop
}$$

There is also an observation to be made. The fan $\Sigma_{\ocM_g}$ is put together based on the strata of $\ocM_g$, namely the spaces $\ocM_\Gamma$. The corresponding cones $\sigma_{\ocM_\Gamma}$ are glued together via the same rule used to glue together the space $\ocM_g^\Trop$. Putting these facts together leads to the following:

\begin{theorem}
We have a canonical isomorphism $\Sigma_X \simeq \ocM_g^\Trop$, making the following diagram commutative: 
$$\xymatrix{ 
 \ocM_g^\An \ar^p[r] \ar_{\Trop}[dr] & \Sigma_X\ar@{=}[d] \\
& \ocM_g^\Trop
}$$
\end{theorem}

In particular the map $\Trop$ is proper and continuous, as claimed.

\providecommand{\bysame}{\leavevmode\hbox to3em{\hrulefill}\thinspace}
\providecommand{\MR}{\relax\ifhmode\unskip\space\fi MR }
\providecommand{\MRhref}[2]{%
  \href{http://www.ams.org/mathscinet-getitem?mr=#1}{#2}
}
\providecommand{\href}[2]{#2}

\end{document}